\documentclass[a4paper,12pt,article]{aplimattpl}

\usepackage[T1]{fontenc}
\usepackage[utf8]{inputenc}

\usepackage{amsmath}
\usepackage{amssymb}
\usepackage{graphicx}

\newtheorem{theorem}{Theorem}

\newtheorem{definition}{Definition}

\begin{document}

\pagestyle{aplimat}


\authorsandtitle{PAVL\'IKOV\'A So\v{n}a/SK, and  \v{S}EV\v{C}OVI\v{C} Daniel/SK}{Invertibility of graphs with a unique perfect matching}{1.2}

\abstract{
In this paper we investigate invertibility of graphs with a unique perfect matching, i.e. graphs having a unique 1-factor. We recall the new notion of the so-called negatively invertible graphs investigated by the authors in the recent paper. It is an extension of the classical definition of an inverse graph due to Godsil. We characterize all graphs with a unique perfect matching on $m\le 6$ vertices with respect to their positive and negative invertibility. We show that negatively invertible graphs exhibit properties like selfinvertibility which cannot be observed within the class of positively invertible non-bipartite graphs with a unique perfect matching.}



\newsavebox{\authors}
\savebox{\authors}{%
\parbox{1.0\textwidth}{%
		\setlength{\parskip}{6pt}%
    \textbf{Pavl\'{\i}kov\'a So\v{n}a, RNDr., CSc.}\\
    Institute of Information Engineering, Automation, and Mathematics, FCFT, Slovak Technical University, 812 37 Bratislava, Slovak republic\\
    sona.pavlikova@stuba.sk
    
    \textbf{\v{S}ev\v{c}ovi\v{c} Daniel, Prof., RNDr., CSc.}\\
    Institute of Mathematics and Physics, Faculty of Mechanical Engineering, Slovak University of Technology in Bratislava, Nám. slobody 17, 812 31 Bratislava, Slovak republic\\
    sevcovic@fmph.uniba.sk    
    }
}

\begin{aplart}

\section{Introduction and preliminaries}

The main purpose of this paper is to investigate invertibility of graphs with a unique perfect matching (graphs having a unique 1-factor) (c.f. Pavl\'\i kov\'a \emph{et al.} \cite{Pavlikova1990,Pavlikova1994}). Following our recent paper \cite{Pavlikova2016}, we recall the new notion of the so-called negatively invertible graphs. It is an extension of the classical definition of an inverse graph due to Godsil \cite{Godsil1985}. Our goal is to characterize all graphs with a unique perfect matching on $m\le 6$ vertices with respect to their positive and negative invertibility.  

According to the result due to Sachs (c.f. \cite[p. 32]{Cvetkovic1988}), it is known that, for a bipartite graph $G$ having a unique perfect matching, the adjacency matrix $A$ is integrally invertible. The aim of this paper is to investigate a larger class of graphs which are not necessarily bipartite but still have a unique 1-factor and an integrally invertible adjacency matrix. 

Let $G=(V,E)$ be an undirected graph with the set of vertices $V$ and the set of edges $E$. By $A_G$ we denote its adjacency matrix. Conversely, if $A$ is a $\{0,1\}$-symmetric matrix then $G_A$ denotes the graph with the adjacency matrix $A$. The spectrum $\sigma(G)$ of an undirected graph $G$ consists of eigenvalues of its adjacency symmetric matrix $A_G$, i.e. $\sigma(G)=\sigma(A_G)=\{\lambda, \lambda$ is an eigenvalue of $A_G\}$ (c.f. \cite{Cvetkovic1978, Cvetkovic1988}). If the spectrum does not contain zero then the adjacency matrix $A$ is invertible. 

Clearly, a graph $G_A$ has the integral inverse matrix $A^{-1}$ if and only if the determinant $\det(A)=\pm1$ (c.f. Kirkland and Akbari \cite{KirklandAkb2007}). The inverse matrix $A^{-1}$ however may contain positive as well as negative entries. In fact, it is well known that the inverse matrix of the adjacency matrix of a graph is nonnegative if and only if the graph is a union of isolated edges (c.f. Harary \cite{Har}). Nevertheless, the concept of an inverse graph is indeed based on the inverse of its adjacency matrix. In  \cite{Godsil1985} Godsil defined a graph to be invertible if the inverse of its (non-singular) adjacency matrix is diagonally similar (c.f. \cite{Zas}) to a non-negative matrix.

\begin{definition}[Godsil \cite{Godsil1985}]
An undirected vertex-labeled graph $G_A$ is called positively invertible (or invertible in Godsil's sense) if inverse of its adjacency matrix is signable to a nonnegative integer matrix, i.e. there exists a diagonal matrix $D^A$ containing $\pm1$ elements and such that all elements of the matrix $D^A A^{-1} D^A$ are nonnegative integers only. The inverse graph $H=G_A^{-1}$ is defined by the adjacency matrix $A_H=D^A A^{-1} D^A$.
\end{definition}

Unfortunately, the class of positively invertible graphs does not contain many important examples of graphs having integral inverse matrix. Among them there is a graph representing the chemical organic molecule of fulvene (see Fig.~\ref{fig-fulvene}). The inverse of its adjacency matrix is integral but it is not signable to a nonnegative matrix.

\begin{figure}
\centering
\includegraphics[width=0.2\textwidth]{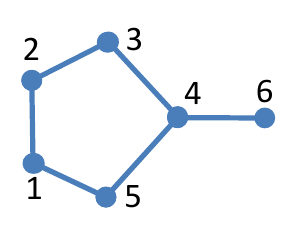}

\caption{A graph representing the structural graph of a chemical molecule of fulvene}
\label{fig-fulvene}

\end{figure}

The structural graph of fulvene however plays an important role in computational chemistry as it represents a chemical molecule. For instance, spectral properties of inverse graphs are important in order to compute the so-called binding energies of molecular orbitals. 

Motivated by this shortcoming of positive invertibility of a graph, in the recent paper \cite{Pavlikova2016}, Pavl\'\i kov\'a and \v{S}ev\v{c}ovi\v{c}  introduced the new concept of the so-called negative invertibility of a graph with integral inverse of its adjacency matrix. 

\begin{definition}[Pavl\'\i kov\'a and \v{S}ev\v{c}ovi\v{c}  \cite{Pavlikova2016}]
An undirected vertex-labeled graph $G_A$ is called negative invertible if the inverse matrix of its adjacency matrix is signable to a nonpositive integer matrix, i.e. there exists a diagonal matrix $D^A$ containing $\pm1$ elements and such that all elements of the matrix $D^A A^{-1} D^A$ are nonpositive integers only. The inverse graph $H=G_A^{-1}$ is defined by the adjacency matrix $A_H=-D^A A^{-1} D^A$.
\end{definition}

The adjacency matrix $A$ of the fulvene graph and its inverse matrix $A^{-1}$ are as follows:
\[
A=
\begin{pmatrix}
0 & 1 & 0 & 0 & 1 & 0 \\ 
1 & 0 & 1 & 0 & 0 & 0 \\ 
0 & 1 & 0 & 1 & 0 & 0 \\ 
0 & 0 & 1 & 0 & 1 & 1 \\ 
1 & 0 & 0 & 1 & 0 & 0 \\ 
0 & 0 & 0 & 1 & 0 & 0
\end{pmatrix},
\qquad 
A^{-1}=
\begin{pmatrix}
0 & 0 & 0 & 0 & 1 & -1 \\ 
0 & 0 & 1 & 0 & 0 & -1 \\ 
0 & 1 & 0 & 0 & -1 & 1 \\ 
0 & 0 & 0 & 0 & 0 & 1 \\ 
1 & 0 & -1 & 0 & 0 & 1 \\ 
-1 & -1 & 1 & 1 & 1 & -2
\end{pmatrix}.
\]
The matrix $A^{-1}$ is signable to a nonpositive matrix by the diagonal matrix $D^A=diag(1,1,-1,-1,-1,1)$.

Both positive as well as negatively invertible graphs have the following useful property. If $G=G_A$ is positively invertible then $A_H=D^A A^{-1} D^A$ is the adjacency matrix of the inverse graph $H=G^{-1}$. As $A= D^A A_H^{-1} D^A$ we have the identity $G=(G^{-1})^{-1}$. This property also remains true in the case of negatively invertible graphs. Indeed, if $G=G_A$ is negatively invertible then $A_H=-D^A A^{-1} D^A$ is the adjacency matrix of the inverse graph $H=G^{-1}$. As $A= - D^A A_H^{-1} D^A$ we again obtain  $G=(G^{-1})^{-1}$.

In the next section we recall the complete list of graphs with a unique perfect matching on $m\le 6$ vertices discovered by  Pavl\'\i kov\'a and \v{S}ev\v{c}ovi\v{c}  \cite{Pavlikova2016}. In addition, we present their inverse graphs and maximal subgraphs with a unique perfect matching. Interestingly enough, there are no positively invertible non-bipartite graphs with a unique perfect matching that are selfinvertible. On the other hand, there exist two negatively invertible graphs which are selfinvertible.

\section{Connected graphs with a unique 1-factor}

In this section we present the complete list of invertible graphs on $m\le 6$ vertices with a unique 1-factor. Recall that a graph $G$ has a unique 1-factor if there exists a unique subgraph spanning $G$ and such that each vertex has the degree $1$. Notice that any graph having a 1-factor should have even number of vertices. 

For $m=2$ the graph $K_2$ is the unique connected graph and it has a unique 1-factor. It is a positively and negatively invertible bipartite graph. 

For $m=4$ there are two connected graphs $Q_1, Q_2$ with a unique 1-factor shown in  Fig.~\ref{fig-quarticfamily}. The graph $Q_1$ is bipartite and selfinvertible. The graph $Q_2$ is positively invertible but it is not selfinvertible.

\begin{figure}[h]

\centering

\ \hskip -1truecm\includegraphics[width=0.16\textwidth]{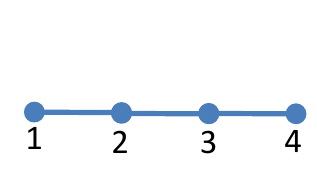}
\hskip1truecm
\includegraphics[width=0.16\textwidth]{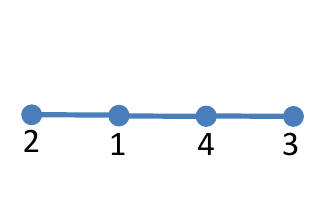}
\hskip1truecm
\includegraphics[width=0.16\textwidth]{figures/quartic-1}
\\
\ \hskip0.5truecm $Q_{1}$ \hskip 2.5truecm $(Q_{1})^{-1}$ \hskip 2truecm $Q_{1} = (Q_{1})^{-1}$

\bigskip

\ \hskip -1truecm\includegraphics[width=0.14\textwidth]{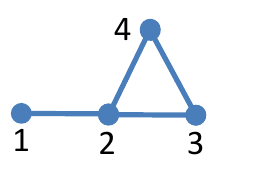}
\hskip1truecm
\includegraphics[width=0.16\textwidth]{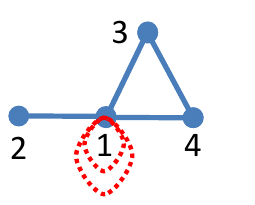}
\hskip1truecm
\includegraphics[width=0.14\textwidth]{figures/quartic-2}
\\
\ \hskip0.5truecm $Q_{2}$ \hskip 2.5truecm $(Q_{2})^{-1}$ \hskip 1.5truecm $Q_{2} \subseteq (Q_{2})^{-1}$

\bigskip

\caption{Invertible graphs on four vertices with a unique perfect matching, their inverse graphs and maximal subgraphs with a unique perfect matching}

\label{fig-quarticfamily}
\end{figure}

Next we analyse the connected graphs on $m=6$ vertices with a unique $1$-factor. We have the following theorem.

\begin{theorem}\cite{Pavlikova2016}\label{theo-hexa}
There exist 20 undirected connected graphs on $m=6$ vertices with a unique 1-factor shown in Figures~\ref{fig-bipartite},\ref{fig-positive1},\ref{fig-positive2},\ref{fig-negative},\ref{fig-noninvertible}. All of them have invertible adjacency matrices. Except of the graph $H_{19}$ they are integrally invertible.

There are three bipartite graphs $H_1, H_2$, and $H_6$ which are simultaneously positively and negatively invertible. There are twelve graphs 
\[
H_3, H_4, H_7, H_8, H_9, H_{13}, H_{14}, H_{15}, H_{16}, H_{17}, H_{18}, H_{20},
\] 
which are positively invertible. The three graphs $H_5, H_{10}$, and $H_{12}$ are negatively invertible. The integrally invertible graph  $H_{11}$ is neither positively nor negatively invertible. The graphs $H_8$, and $H_{18}$ are iso-spectral but not isomorphic. 

\end{theorem}

Recall that the proof of this result is based on the well-known Kotzig's theorem. It states that a graph with a unique 1-factor should contain a bridge belonging to the 1-factor sub-graph. Using this characterization one can construct all 20 connected graphs on 6 vertices having a unique 1-factor. We refer the reader to \cite{Pavlikova2016} for details.

\begin{table}
\small
\caption{The family of graphs on 6 vertices with a unique 1-factor and their signability. Source: \cite{Pavlikova2016} } 
\begin{center}
\begin{tabular}{c||c}
Graph & invertibility  \\
\hline\hline
   $H_1$ &  positively and negatively invertible (bipartite) \\
\hline
   $H_2$ &  positively and negatively invertible (bipartite) \\
\hline
   $H_3$ &  positively  invertible  \\ 
\hline
   $H_4$ &  positively  invertible  \\
\hline
   $H_5$ &  negatively invertible \\
\hline
   $H_6$ &  positively and negatively invertible (bipartite) \\ 
\hline
   $H_7$ &  positively  invertible  \\
\hline
   $H_8$ &  positively  invertible  \\
\hline
   $H_9$ &  positively  invertible  \\
\hline
$H_{10}$ &  negatively invertible \\
\hline
$H_{11}$ &  noninvertible with integral inverse of the adjacency matrix\\ 
\hline
$H_{12}$ &  negatively invertible \\
\hline
$H_{13}$ &  positively  invertible  \\
\hline
$H_{14}$ &  positively  invertible  \\
\hline
$H_{15}$ &  positively  invertible  \\
\hline
$H_{16}$ &  positively  invertible  \\
\hline
$H_{17}$ &  positively  invertible  \\
\hline
$H_{18}$ &  positively  invertible  \\
\hline
$H_{19}$ &  noninvertible with nonintegral inverse of the adjacency matrix\\
\hline
$H_{20}$ &  positively  invertible  \\
\hline
\hline\hline
\end{tabular}

\end{center}
\end{table}

\subsection{Positively and negatively invertible (bipartite) graphs with a unique perfect matching on $m=6$ vertices}

Recall that a graph $G_B$ is called bipartite if the set of vertices can be divided into two disjoint subsets such that every edge connects a vertex from the first subset into the one from the second subset.

\begin{theorem}[Pavl\'\i kov\'a and \v{S}ev\v{c}ovi\v{c} \cite{Pavlikova2016}]
\label{theo-bipartite}
Let $G$ be an integrally invertible graph. Then  $G$ is a bipartite graph if and only if $G$ is  simultaneously positively and negatively invertible.
\end{theorem}

In \cite{McLeman2014} McLeman and McNicholas considered graphs which  can be obtained from a given connected graph by attaching pendant edges to each of vertices. They proved that such graphs are selfinvertible, i.e. $H^{-1}=H$. An example illustrating their result is the graph $H_2$ which is obtained from a path with three vertices by adding pendant edges to each of them.

\begin{figure}[h!]

\centering

\ \hskip -2.5truecm\includegraphics[width=0.3\textwidth]{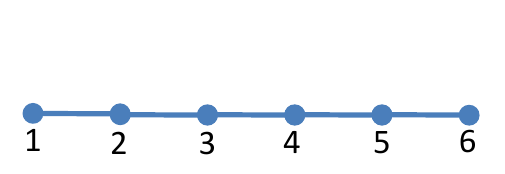}
\hskip1truecm
\includegraphics[width=0.14\textwidth]{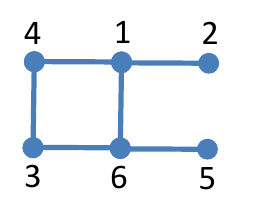}
\hskip1truecm
\includegraphics[width=0.12\textwidth]{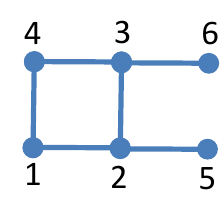}
\\
\ \hskip0truecm $H_1$ \hskip 4truecm $(H_1)^{-1}$ \hskip 1.5truecm $H_6 = (H_1)^{-1}$

\bigskip

\ \hskip -1truecm\includegraphics[width=0.22\textwidth]{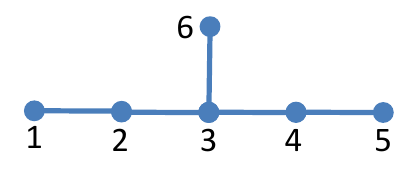}
\hskip1truecm
\includegraphics[width=0.22\textwidth]{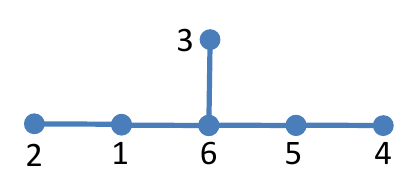}
\hskip1truecm
\includegraphics[width=0.22\textwidth]{figures/hexa-02}
\\
\ \hskip0truecm $H_2$ \hskip 3.5truecm $(H_2)^{-1}$ \hskip 2.5truecm $H_2 = (H_2)^{-1}$

\bigskip

\ \hskip -0.7truecm\includegraphics[width=0.12\textwidth]{figures/hexa-06}
\hskip1truecm
\includegraphics[width=0.26\textwidth]{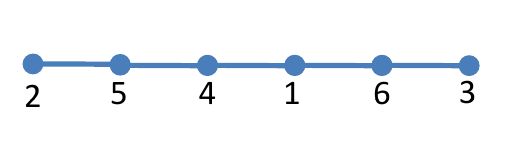}
\hskip1truecm
\includegraphics[width=0.26\textwidth]{figures/hexa-01}
\\
\ \hskip-0.5truecm $H_6$ \hskip 3truecm $(H_6)^{-1}$ \hskip 3truecm $H_1 = (H_6)^{-1}$

\caption{Positively and negatively invertible (bipartite) graphs  with a unique perfect matching, their inverses and maximal subgraphs with a unique perfect matching}
\label{fig-bipartite}
\end{figure}

In  Fig.~\ref{fig-bipartite} we show all three bipartite graphs on $m=6$ vertices  with a unique perfect matching. According to the previous theorem they are simultaneously positive as well as negative invertible.

\subsection{Positively invertible graphs with a unique perfect matching on $m=6$ vertices}

In this section we present the complete list of twelve non-bipartite graphs with a unique perfect matching on $m=6$ vertices which are positively invariant (see Fig.~\ref{fig-positive1},\ref{fig-positive2}). Interestingly, none of them is selfinvertible and their inverse graphs contain vertices with loops with even multiplicity. The inverse  graphs $(H_{4})^{-1},(H_{8})^{-1},(H_{13})^{-1}, (H_{15})^{-1}$ contain multiple edges. There are four graphs denoted by $H_9, H_{13}, H_{14}$, and $H_{18}$ such that they are maximal subgraphs of their inverse graphs, i.e. $H_i\subseteq (H_i)^{-1}, i=9,13,14,18$. There are two pairs of graphs $(H_i, H_j)$ such that $H_i\subseteq (H_j)^{-1}$ and $H_j\subseteq (H_i)^{-1}$ for pairs $(i,j)\in \{ (3,7), (15,17)\}$.

\begin{figure}[h!]

\centering

\ \hskip -1truecm\includegraphics[width=0.22\textwidth]{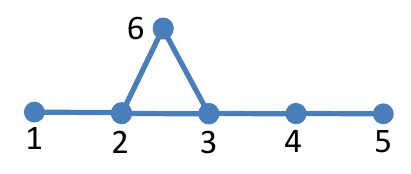}
\hskip1truecm
\includegraphics[width=0.12\textwidth]{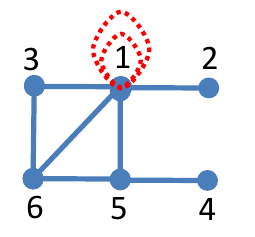}
\hskip1truecm
\includegraphics[width=0.12\textwidth]{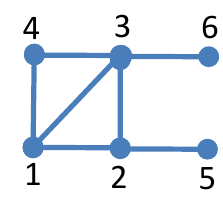}
\\
\ \hskip1truecm $H_3$ \hskip 2.5truecm $(H_3)^{-1}$ \hskip 1.5truecm $H_7 \subseteq (H_3)^{-1}$

\bigskip

\ \hskip -1truecm\includegraphics[width=0.22\textwidth]{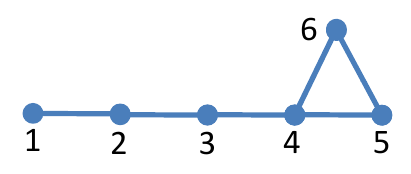}
\hskip1truecm
\includegraphics[width=0.12\textwidth]{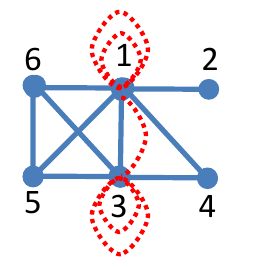}
\hskip1truecm
\includegraphics[width=0.12\textwidth]{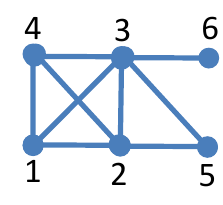}
\\
\ \hskip1.5truecm $H_4$ \hskip 2.5truecm $(H_4)^{-1}$ \hskip 1.5truecm $H_{17} \subseteq (H_{4})^{-1}$

\bigskip

\ \hskip -1truecm\includegraphics[width=0.11\textwidth]{figures/hexa-07}
\hskip1truecm
\includegraphics[width=0.22\textwidth]{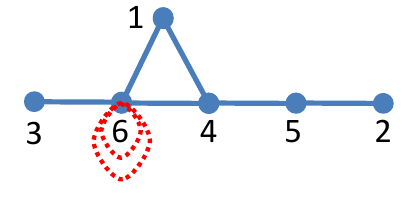}
\hskip1truecm
\includegraphics[width=0.22\textwidth]{figures/hexa-03}
\\
\ \hskip-0.5truecm $H_7$ \hskip 2.5truecm $(H_7)^{-1}$ \hskip 3truecm $H_{3} \subseteq (H_{7})^{-1}$

\bigskip

\ \hskip -1truecm\includegraphics[width=0.16\textwidth]{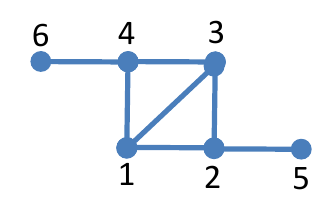}
\hskip1truecm
\includegraphics[width=0.13\textwidth]{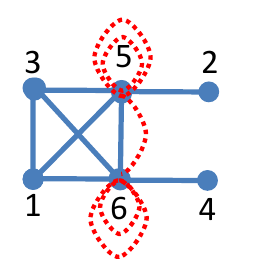}
\hskip1truecm
\includegraphics[width=0.12\textwidth]{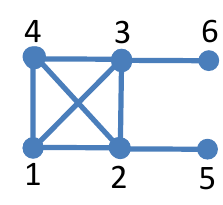}
\\
\ \hskip1truecm $H_8$ \hskip 2truecm $(H_8)^{-1}$ \hskip 1.5truecm $H_{9} \subseteq (H_{8})^{-1}$

\bigskip

\ \hskip -1truecm\includegraphics[width=0.11\textwidth]{figures/hexa-09}
\hskip1truecm
\includegraphics[width=0.13\textwidth]{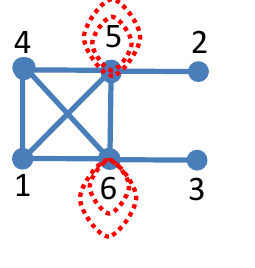}
\hskip1truecm
\includegraphics[width=0.11\textwidth]{figures/hexa-09}
\\
\ \hskip0.5truecm $H_9$ \hskip 2truecm $(H_9)^{-1}$ \hskip 1.5truecm $H_{9} \subseteq (H_{9})^{-1}$

\bigskip

\ \hskip -1truecm\includegraphics[width=0.14\textwidth]{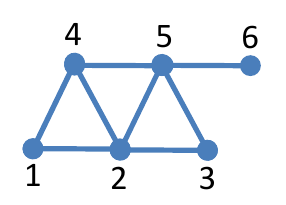}
\hskip1truecm
\includegraphics[width=0.14\textwidth]{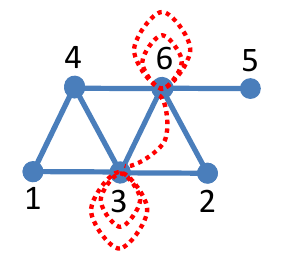}
\hskip1truecm
\includegraphics[width=0.14\textwidth]{figures/hexa-13}
\\
\ \hskip0truecm $H_{13}$ \hskip 2.5truecm $(H_{13})^{-1}$ \hskip 1.5truecm $H_{13} \subseteq (H_{13})^{-1}$

\caption{Positively invertible graphs  with a unique perfect matching, their inverses and maximal subgraphs with a unique perfect matching (1. part)}
\label{fig-positive1}

\end{figure}

\begin{figure}[h!]

\centering

\ \hskip -1truecm\includegraphics[width=0.11\textwidth]{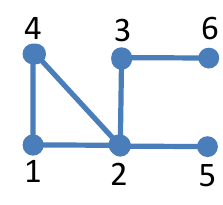}
\hskip1truecm
\includegraphics[width=0.13\textwidth]{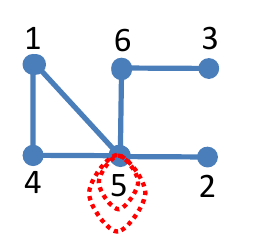}
\hskip1truecm
\includegraphics[width=0.12\textwidth]{figures/hexa-14}
\\
\ \hskip0.5truecm $H_{14}$ \hskip 2truecm $(H_{14})^{-1}$ \hskip 1truecm $H_{14} \subseteq (H_{14})^{-1}$

\bigskip

\ \hskip -1truecm\includegraphics[width=0.14\textwidth]{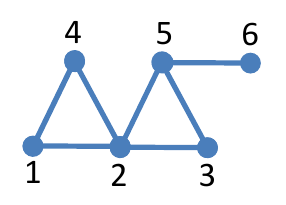}
\hskip1truecm
\includegraphics[width=0.13\textwidth]{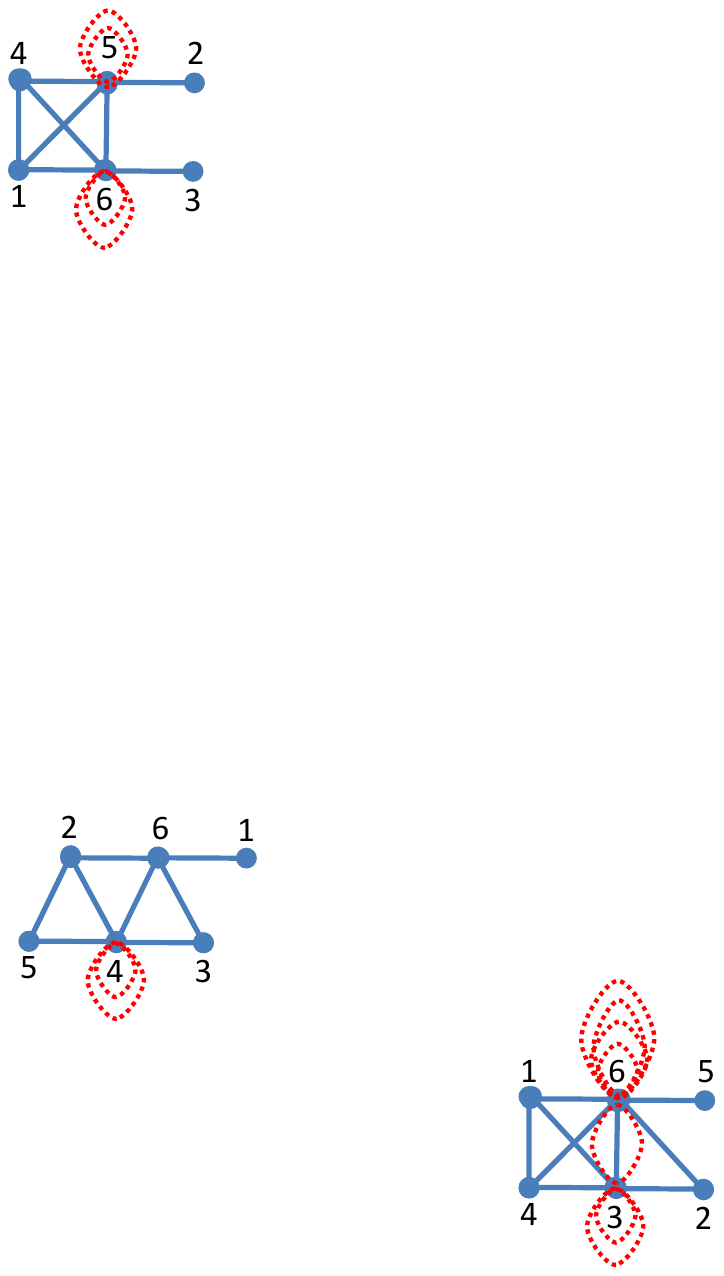}
\hskip1truecm
\includegraphics[width=0.12\textwidth]{figures/hexa-17}
\\
\ \hskip0.5truecm $H_{15}$ \hskip 2truecm $(H_{15})^{-1}$ \hskip 1truecm $H_{17} \subseteq (H_{15})^{-1}$

\bigskip

\ \hskip -1truecm\includegraphics[width=0.15\textwidth]{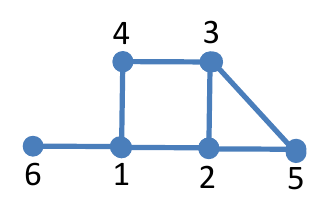}
\hskip1truecm
\includegraphics[width=0.14\textwidth]{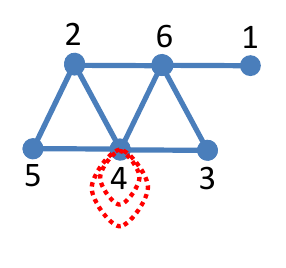}
\hskip1truecm
\includegraphics[width=0.14\textwidth]{figures/hexa-13}
\\
\ \hskip0.5truecm $H_{16}$ \hskip 2truecm $(H_{16})^{-1}$ \hskip 1truecm $H_{13} \subseteq (H_{16})^{-1}$

\bigskip

\ \hskip -1truecm\includegraphics[width=0.12\textwidth]{figures/hexa-17}
\hskip1truecm
\includegraphics[width=0.15\textwidth]{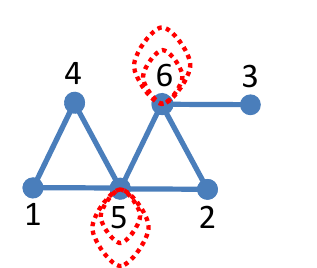}
\hskip1truecm
\includegraphics[width=0.14\textwidth]{figures/hexa-15}
\\
\ \hskip0.5truecm $H_{17}$ \hskip 2truecm $(H_{17})^{-1}$ \hskip 1truecm $H_{15} \subseteq (H_{17})^{-1}$

\bigskip

\ \hskip -1truecm\includegraphics[width=0.12\textwidth]{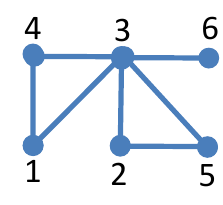}
\hskip1truecm
\includegraphics[width=0.13\textwidth]{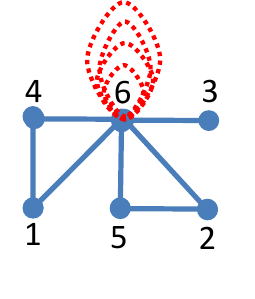}
\hskip1truecm
\includegraphics[width=0.11\textwidth]{figures/hexa-18}
\\
\ \hskip0.5truecm $H_{18}$ \hskip 2truecm $(H_{18})^{-1}$ \hskip 1truecm $H_{18} \subseteq (H_{18})^{-1}$

\bigskip

\ \hskip -1truecm\includegraphics[width=0.12\textwidth]{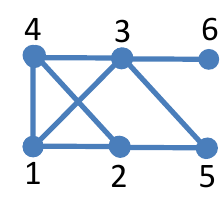}
\hskip1truecm
\includegraphics[width=0.23\textwidth]{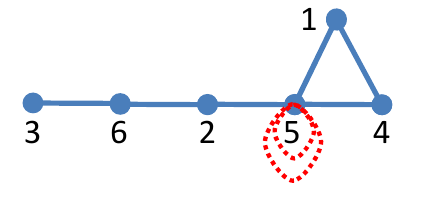}
\hskip1truecm
\includegraphics[width=0.23\textwidth]{figures/hexa-04}
\\
\ \hskip-0.5truecm $H_{20}$ \hskip 2.5truecm $(H_{20})^{-1}$ \hskip 3truecm $H_{4} \subseteq (H_{20})^{-1}$

\caption{Positively invertible graphs  with a unique perfect matching, their inverses and maximal subgraphs with a unique perfect matching (2. part)}
\label{fig-positive2}

\end{figure}

\subsection{Negatively invertible graphs with a unique perfect matching on $m=6$ vertices}

In this section we present the complete list of all three non-bipartite graphs with a unique perfect matching on $m=6$ vertices which are negatively invariant. The graphs $H_5$, and $H_{12}$ are selfinvertible. Recall that according to the results due to  McLeman and McNicholas \cite{McLeman2014} a graph which is obtained from a given connected graph by attaching pendant edges to each of vertices is selfinvertible. Among negatively invertible graphs with unique perfect matching on $m=6$ vertices as an example illustrating this property one can present the graph $H_5$. In contrast to positively invertible bipartite graphs there exists the graph $H_{12}$ which is selfinvertible but it cannot be obtained by means of the construction due to McLeman and McNicholas.

\begin{figure}[h]

\centering

\ \hskip -0.5truecm\includegraphics[width=0.14\textwidth]{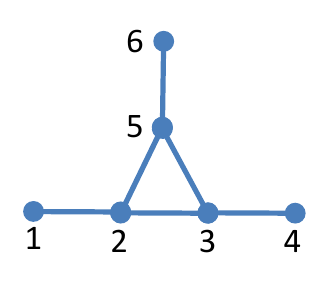}
\hskip1truecm
\includegraphics[width=0.16\textwidth]{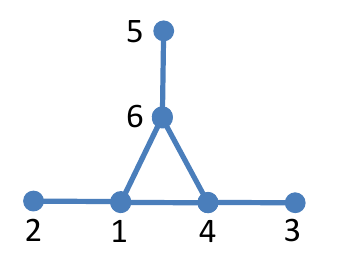}
\hskip1truecm
\includegraphics[width=0.14\textwidth]{figures/hexa-05}
\\
\ \hskip0.5truecm $H_{5}$ \hskip2.5truecm $(H_{5})^{-1}$ \hskip 2truecm $H_{5} = (H_{5})^{-1}$

\bigskip

\ \hskip -1truecm\includegraphics[width=0.14\textwidth]{figures/hexa-10}
\hskip1truecm
\includegraphics[width=0.16\textwidth]{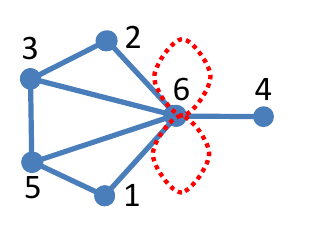}
\hskip1truecm
\includegraphics[width=0.14\textwidth]{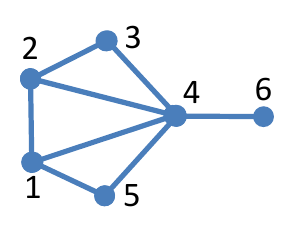}
\\
\ \hskip0.5truecm $H_{10}$ \hskip 2.5truecm $(H_{10})^{-1}$ \hskip 2truecm $H_{11} \subseteq (H_{10})^{-1}$

\bigskip

\ \hskip -1truecm\includegraphics[width=0.14\textwidth]{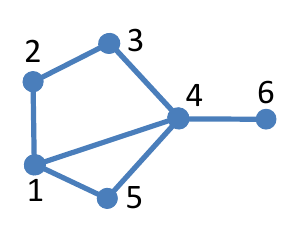}
\hskip1truecm
\includegraphics[width=0.16\textwidth]{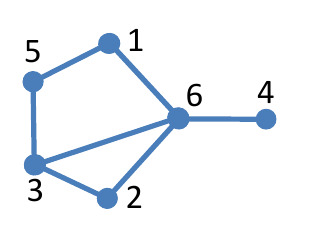}
\hskip1truecm
\includegraphics[width=0.14\textwidth]{figures/hexa-12}
\\
\ \hskip0.5truecm $H_{12}$ \hskip 2.5truecm $(H_{12})^{-1}$ \hskip 2truecm $H_{12} = (H_{12})^{-1}$

\caption{Negatively invertible graphs with a unique perfect matching, their inverses and maximal subgraphs with a unique perfect matching}
\label{fig-negative}

\end{figure}

\subsection{Noninvertible graphs with a unique perfect matching on $m=6$ vertices}

In this section we present two remaining graphs on $m=6$ vertices which are not invertible. The first graph is denoted by $H_{11}$ and it has the adjacency matrix which is integrally invertible. Nevertheless, $A^{-1}$ is neither positively nor negatively signable, so $H_{11}$ is not invertible. The graph $H_{19}$ has the adjacency matrix with $\det(A)=3$ and, as a consequence, $A^{-1}$ is not integer matrix. Hence $H_{19}$ is noninvertible as well.

\begin{figure}[ht]

\centering

\ \hskip -1truecm\includegraphics[width=0.14\textwidth]{figures/hexa-11}
\hskip1truecm
\includegraphics[width=0.17\textwidth]{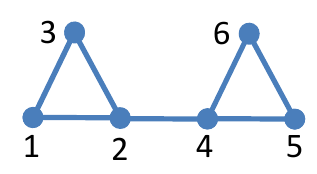}
\qquad\quad
\\
\  \hskip -1truecm $H_{11}$ \hskip 3truecm $H_{19}$

\caption{Noninvertible graphs  with a unique perfect matching: $H_{11}$ has integral inverse of its adjacency matrix but it is neither positively nor negatively invertible, $H_{19}$ is not integrally invertible}
\label{fig-noninvertible}

\end{figure}

\section{Conclusions}
In this note we investigated invertible graphs with a unique perfect matching on $m\le 6$ vertices. We recalled classical concept of graph inversion due to Godsil which we referred to as positive invertibility of a graph. We also recalled the new concept of negative invertibility proposed by the authors in \cite{Pavlikova2016}. By inspecting properties of the inverse graphs we showed that negatively invertible graphs exhibits properties like selfinvertibility which cannot be observed within the class of positively invertible non-bipartite graphs with a unique perfect matching on $m\le 6$ vertices.

\section*{Acknowledgements} The research was supported by the APVV Research Grant 0136-12 (SP) and VEGA grant 1/0780/15 (D\v{S}).

\section*{Current address}

\vspace{1.2ex}
\usebox{\authors}

\end{aplart}
\end{document}